\documentclass[12pt]{amsart}
\usepackage{amsmath,amsthm}

\def\transpose#1{\kern1pt{^t\kern-3pt#1}}%
\def\Fix{{\rm Fix}}
\def\pmatrix{\left(\begin{matrix}}
\def\endpmatrix{\end{matrix}\right)}

\def\Sp{\operatorname{Sp}}
\def\id{{\rm id}}

\def\Z{{\mathbb Z}}

\def\C{{\mathbb C}}
\def\N{{\mathbb N}}

\def\Pj{{\mathbb P}}

\def\cal{\mathcal}

\def\<{\langle}
\def\>{\rangle}

\def\t{\theta}

\def\b{\beta}

\def\H{{\mathbb H}}

\def\tch#1#2{{\left[\begin{matrix}#1\\ #2\end{matrix}\right]}}
\def\tt#1#2{{\t\tch{#1}{#2}}}

\theoremstyle{plain}
\newtheorem{thm}{Theorem}
\newtheorem{lm}[thm]{Lemma}
\newtheorem{prop}[thm]{Proposition}

\newtheorem{rem}[thm]{Remark}
\theoremstyle{definition}
\begin{document}
\title[Calabi-Yau  Siegel threefold]{The geometry and arithmetic of a Calabi-
Yau  Siegel threefold}
\author[Cynk, Freitag \and Salvati Manni]{S. Cynk,
 E. Freitag \and R. Salvati Manni}
\address{Institute of Mathematics\\Jagiellonian University      
ul.  \L ojasiewicza  6\\30-348 Krak\'ow\\Polska} 
\email{slawomir.cynk@uj.edu.pl}
\address{Mathematisches Institut\\Im Neuenheimer Feld 288\\D69120 Heidelberg}
\email{freitag@mathi.uni-heidelberg.de}
\address{Dipartimento di Matematica\\Universit\`a La Sapienza\\Piazzale Aldo 
Moro, 2- I 00185 Roma}
\email{salvati@mat.uniroma1.it}
\date{\today}
\maketitle

\section{Introduction}

\noindent
In two recent papers  \cite{FS1} and  \cite{FS2}, the last two named
authors described Siegel modular varieties  which admit a Calabi-Yau
model. They used    two different methods, but essentially they
restrict to consider the action of a  finite group  G, fixing a
holomorphic three form, on a smooth projective  variety  $M$.  In the
first case the variety $M$ was a toroidal compactification of the
Siegel modular variety of level 4, in the second case they started
from a small resolution of  a singular  Siegel modular variety $\cal
X$ introduced by van Geemen and Nygaard, cf. \cite{vGN} . The second
method appears more powerful and leads to the construction of  more
than 4000 Calabi Yau varieties of which one can compute  Hodge
numbers.  
 This will be the content of a forthcoming paper. \smallskip

A careful analysis of the first method leads to introduce a modular
variety of particular interest related to a significant   modular
group. The aim of this paper is to treat in details the associated
modular  variety $\cal Y$ that  has a Calabi-Yau model, $\tilde{\cal
  Y}$. We shall describe  its geometry and the structure of the ring
of modular forms using several approaches. We shall illustrate two
different methods of producing the Hodge numbers. The first uses the
definition of  $\cal Y$ as the quotient  of 
$\cal X$  modulo a finite group $K$.  In this second case 
we  will get the Hodge numbers  considering the action of the group
$K$ on a crepant  resolution $\tilde{\cal X}$ of $\cal X$.\medskip

The  second, purely
algebraic geometric, uses the  equations derived from the ring of 
modular forms and is based on determining explicitly the
Calabi-Yau model $\tilde{\cal Y}$ and computing the Picard group and
the Euler characteristic.\smallskip

\section{ Modular varieties}

\noindent
As in   \cite{FS2} the staring point of our investigation is  the variety 
\[\mathcal X:\qquad\qquad\qquad\begin{array}{c}
Y_{0}^{2}=X_{0}^{2}+X_{1}^{2}+X_{2}^{2}+X_{3}^{2}\\
 Y_{1}^{2}=X_{0}^{2}-X_{1}^{2}+X_{2}^{2}-X_{3}^{2}\\
 Y_{2}^{2}=X_{0}^{2}+X_{1}^{2}-X_{2}^{2}-X_{3}^{2}\\
Y_{3}^{2}=X_{0}^{2}-X_{1}^{2}-X_{2}^{2}+X_{3}^{2}\\
\end{array}\qquad\qquad
\]
This is a modular  variety , in the sense that is biholomorphic to 
the Satake compactification of $\H_2/\Gamma'$ for a 
certain subgroup $\Gamma' \subset\Sp(4, \Z)$.
For details, we refer to  \cite{vGN}, \cite{CM} and \cite{FS2}, we just recall the basic informations that we need.\smallskip
  
Let  $\H_n$ be  the Siegel
upper half space of symmetric complex matrices with positive definite imaginary part. The symplectic group $\Gamma_n:=\Sp(2n,\Z)$ acts on $\H_n$  via
$$
  \pmatrix A&B\\  C&D\endpmatrix\cdot Z:= (AZ+B)(CZ+D)^{-1} .
$$
Here we think of elements of $\Gamma_n$ as consisting of four
$n\times n$ blocks. 
For any subgroup  of finite index $\Gamma\subset \Gamma_n$  
the Satake compactification  $\overline{\H_n/\Gamma}$ of  the quotient $\H_n/\Gamma$  
 is the  projective variety associated
to a graded algebra of modular forms. 
We recall briefly its definition. A modular form $f$ of weight $r/2$, $r\in\N$,
is a holomorphic function $f$ on $\H_n$ with the transformation property
$$f(MZ)=v(M)\sqrt{\det(CZ+D)}^rf(Z)$$
for all $M\in\Gamma$.
In the case $n=1$ a regularity condition at the cusps has to be added.
Here $v(M)$ is a multiplier system. Essentially it  fulfills the  cocycle condition.
We denote this space by
$[\Gamma,r/2,v]$.
Fixing some starting weight $r_0$ and a multiplier system $v$ for it, we define the
ring
$$A(\Gamma):=\bigoplus_{r\in\N}[\Gamma,rr_0/2,v^r].$$
This turns out to be a finitely generated graded algebra and its associated
projective variety Proj$(A(\Gamma))$  can be identified with the Satake compactification.
The ring depends on the starting weight and the multiplier system but the associated
projective variety does not.
\smallskip

The simplest  examples of modular forms are given by theta constants.  
A characteristic is an element $m={a\choose b}$ from
$(\Z/2\Z)^{2n}$. Here $a,b\in(\Z/2\Z)^n$
are column vectors. The characteristic is called even
if $^tab=0$ and odd otherwise.
The group $\Sp(2n,\Z/2\Z)$ acts on the set of characteristics by
$$M\{m\}:={\transpose M}^{-1}m+\pmatrix(A{\transpose B})_0\\ (C\transpose  D)_0\endpmatrix.$$
Here $S_0$ denotes the column built of the diagonal of a square matrix $S$.
It is well-known that $\Sp(2n,\Z/2\Z)$ acts transitively on the subsets of
even and odd characteristics.
Recall that to any characteristic the theta function
$$\vartheta[m]=\sum_{g\in\Z^n}e^{i\pi  (Z[g+a/2]+^tb(g+a/2))}
\qquad (Z[g]= \,\transpose gZg)$$
can be defined.
Here we use the identification of $\Z/2\Z$ with the
subset $\{0,1\}\subset\Z$. It vanishes if and only if $m$ is odd.
Recall also that the formula
$$\vartheta[M\{m\}](MZ)=v(M,m)\sqrt{\det (CZ+D)}\vartheta[m](Z)$$
holds for $M\in\Gamma_n$, where $v(M,m)$ is a rather delicate
eighth root of unity which depends on the choice of the
square root.
Sometimes, when $n=2$, we will use the notation
$$\vartheta[m]=\vartheta\Bigl[{a_1a_2\atop b_1b_2}\Bigr]
\quad {\rm for}\quad m=\pmatrix a_1\\ a_2\\ b_1\\b_2\endpmatrix.$$
We consider the 8 functions
\medskip

\halign{\qquad$\displaystyle#$\quad\hfil&$\displaystyle#$\quad\hfil
&$\displaystyle#$\quad\hfil&$\displaystyle#$\quad\hfil\cr
\vartheta\Bigl[{00\atop 00}\Bigr](Z),& \vartheta\Bigl[{00\atop 10}\Bigr](Z),&
\vartheta\Bigl[{00\atop 01}\Bigr](Z),& \vartheta\Bigl[{00\atop 11}\Bigr](Z),\cr
\noalign{\vskip3mm}
\vartheta\Bigl[{00\atop 00}\Bigr](2Z),& \vartheta\Bigl[{10\atop 00}\Bigr](2Z),&
\vartheta\Bigl[{01\atop 00}\Bigr](2Z),& \vartheta\Bigl[{11\atop 00}\Bigr](2Z).\cr}

\medskip\noindent
If  we denote them by $Y_0,\dots, Y_3,X_0\dots,X_3$, then 
classical addition formulas  for theta constants show that 
the relations defining $ \cal X$  hold.
These eight forms are modular
forms of weight $1/2$ for a  group $\Gamma'$  that we are  going to define.
\smallskip

We set 
 \begin{eqnarray*}
 \Gamma_n[q]&=&{\rm kernel}\, ( \Gamma_n\to \Sp(2n,\Z/q\Z)),\\
 \Gamma_n[q,2q]&=&\{ M \in\Gamma_n[q]; \quad
(A\transpose  B/q)_0 \equiv (C\transpose  D/q)_0\equiv 0\,  {\rm mod}\, 2 \}, \\
 \Gamma_{n, 0}[q]&=&\{M\in  \Gamma_n; \, C\equiv 0\, {\rm mod}\, q\},\\
\Gamma_{n,0,\vartheta}[q]&=&\{M\in\Gamma_{n, 0}[q] \quad  
(C\transpose  D/q)_0  \equiv 0\,  {\rm mod}\, 2 \}.\\
\end{eqnarray*}
Here $S_0$ denotes the diagonal of the matrix $S$.

\smallskip\noindent
The group $\Gamma'$, which belongs to van Geemen's and Nygaard's variety
is defined by
 $$\Gamma'=\{M\in \Gamma_2[2,4]\cap\Gamma_{2,0,\vartheta}[4];
\quad \det D\equiv\pm1\;\mod\;8\}.$$

\smallskip\noindent  
We are going to recall the main result of  \cite{FS2}.
The group $\Gamma_{n,0}[q]$ can be extended by the Fricke involution
$$ J_q=\pmatrix 0&E/\sqrt q\\  -\sqrt q E&0\endpmatrix.$$
We denote by $\hat{\Gamma}_{2, 0}[2]$ the extension of $\Gamma_{2, 0}[2]$ by $J_2$, i.e.
$$\hat{\Gamma}_{2, 0}[2] =\Gamma_{2, 0}[2]\cup J_2\Gamma_{2, 0}[2].$$
 $\hat{\Gamma}_{2,0}[2]_{\bf n}$ is a subgroup
of index two of $\hat{\Gamma}_{2,0}[2]$ that is  the kernel of a certain
character $\chi_{\bf n}$ that has been explained in \cite{FS1}.  
With these notations we have:
\begin{thm}\label{thm:cymodel}
The Siegel modular threefold, which belongs to a group between
$\Gamma'$ and $\hat{\Gamma}_{2,0}[2]_{\bf n}$, admits a Calabi-Yau model,
more precisely: There exists a desingularization of the
Satake compactification which is a (projective) Calabi-Yau manifold.
\end{thm}
\noindent
So there are thousands of conjugacy classes of intermediate groups, which all
lead to Calabi-Yau manifolds.\smallskip

\section{The variety  $\mathcal Y$ }

\noindent There is one intermediate group of particular interest, namely the group
$$\Gamma=\Gamma_2[2]\cap\Gamma_{2,0}[4].$$
 This group  contains $\Gamma'$ as subgroup of index 32. It is  stable under the Fricke involution
 $J_2$,  for this group (and as a consequence for all groups between it
 and $\hat {\Gamma}_{2,0}[2]_{\bf n}$) we have a completely different
  proof which rests
 on the paper \cite{FS1} and gives a very explicit description of the
 Calabi-Yau model, namely:
 
\begin{thm}
Let $\tilde X(4)$ be the Igusa desingularization of the Satake compactification
 of $\H_2/\Gamma[4]$. Then the quotient
$\tilde X(4)/(\Gamma_2[2]\cap\Gamma_{2,0}[4])$
admits a desingularization, which is a Calabi-Yau manifold 
 \end{thm}
\noindent
For a proof we proceed as it follows. We have to consider translation matrices
$$T_S=\pmatrix E&S\\  0&E\endpmatrix $$
of level two, $S\equiv 0$ mod 2. Such a translation matrix is called
{\it reflective\/} if $S$ is congruent 0 mod 4 to one of the three
$$\pmatrix 2&0\\0&0\endpmatrix ,\quad \pmatrix 0&0\\0&2\endpmatrix ,\quad \pmatrix 2&2\\2&2\endpmatrix .$$
Actually reflective translations act as reflections on the Igusa desingularization
of level four.

\begin{lm}
 The group
$\Gamma $ is generated by
\vskip1mm
\item{\rm 1)} The group $\Gamma_2[4]$,
\item{\rm 2)} The elements of $\hat\Gamma_{2,0}[2]_{\bf n}$, which are
conjugate inside $\Gamma_2$ to the diagonal matrix with diagonal
$(1,-1,1,-1)$.
\item{\rm 3)} All elements of $\hat\Gamma_{2,0}[2]_{\bf n}$, which are
conjugate inside $\Gamma_2$ to a reflective
translation matrix ${E\,S\choose0\,E}$ of
$\Gamma_2[2]$.
\vskip0pt
\end{lm}
\noindent
The proof  can be easily done with the help of a computer.
\smallskip

The lemma  is similar to lemma 1.4 in \cite{FS1}. There the group
$\Gamma_{2,0}[2]_{\bf n}\cap\Gamma_2[2]$ has been characterized by the
same properties 1)--3) with the only difference that the word
``reflective'' has been skipped. The same proof as in [FS] works with this weaker
assumption and gives the result 
that the quotient of the Igusa desingularization for the principal congruence subgroup
of level four
$\tilde X(4)/(\Gamma_2[2]\cap\Gamma_{2,0}[4])$ admits a desingularization
that is
a Calabi-Yau manifold. The same then is true for any group between
$\Gamma $ and $\hat\Gamma_{2,0}[2]_{\bf n}$.\smallskip

\smallskip
By standard method ( going down process)
we can  produce the structure of the ring of modular forms of  this distinguished case:
\begin{prop}
The  ring of modular forms
of even weight  
for the group $\Gamma_0[4]\cap\Gamma[2]$
is generated by
$$ \tt{00}{00} \tt{00}{01}  \tt{00}{10} \tt{00}{11},$$
all pairs  of the form
$$ \tt{00}{ab}^2 \tt{00}{cd}^2$$
and
the  ten even $\theta[ m]^4$.
 
\noindent If one wants the generators also in the  odd weights,
it is enough to add the form   of weight 3
$$T=\tt{10}{00}\tt{10}{01}\tt{01}{00}\tt{01}{10}\tt{11}{00}\tt{11}{11} . $$
\end{prop}
\noindent
To simplify the equations we consider the ring of forms of even weights:
\begin{prop}
{The  ring $A(\Gamma_0[4]\cap\Gamma[2] )^{(2)}$ in the even weights
 is equal to
$$ \C\Bigl[ \tt{00}{00} \tt{00}{01}  \tt{00}{10} \tt{00}{11},  \tt{00}{
00} ^2,  \tt{00}{01}^2,   \tt{00}{10}^2,  \tt{00}{11}^2, y_4 \Bigr]^{(2)}$$
with
$$ y_4=-\tt{10}{
01}^4 -\tt{00}{
11}^4$$}
\end{prop}
\noindent
Denoting the above variables by  $y_5, x_0, x_1, x_2, x_3$ we have the ring
$$ \C[y_5,  x_0, x_1, x_2, x_3 , y_4 ]^{(2)}$$
with  $x_i$ of weight 1 and $y_j$ of weight 2.
We have  also the following defining relations
\begin{equation}
   \label{eq:equations}
\begin{array}{rcl}
  y_{5}^{2}&=&x_{0}x_{1}x_{2}x_{3},\\
  2y_{5}^{2}&=&x_{0}^{2}x_{1}^{2}+x_{0}^{2}x_{3}^{2}+x_{1}^{2}x_{3}^{2}+
  (-x_{2}^{2}+x_{0}^{2}+x_{1}^{2}+x_{3}^{2}+y_{4})y_{4}.
\end{array}
  \end{equation}
We shall denote by $\mathcal Y$ the modular  variety defined by the above equations.\smallskip

We want to explain how we can compute the Hodge  numbers of  a  Calabi--Yau model of the  variety $\mathcal Y$ without the description
of an explicit crepant resolution.

We go back to the modular approach. We need some information about the group
$K:=\Gamma/\Gamma'$.
The basic information is that $K$ is abelian of order 32 and that
all elements are of order two.
So their fixed point loci are known from \cite{FS2}.
We know that they all
extend to a small resolution $\tilde{\mathcal X}$.
We also know  that the fixed point locus is a curve $C\subset \tilde{\mathcal X}$.
The image of $C$ in  $\tilde{\mathcal X}/K$ is the singular locus.
The local structure of a singularity is of the type
$\C^3/A$, where $A$ either is a group of order 2, generated by
a transformation, which changes two signs or the group of order 4 which
contains all sign changes at two positions. It is easy to describe
the crepant resolution for these singularities (see \cite{FS1}) and from
this description on can see:
\begin{lm}The number of exceptional divisors of a crepant resolution of
  $\tilde{\mathcal X}/K$ equals the number of irreducible components of
the fixed point locus of $K$ on  ${\mathcal X}$, modulo $K$.
 \end{lm}
\noindent
One can check that $K$ contains 6 elements which have nodes
as  isolated fixed points. Each of them fixes 16 nodes.
So each node occurs as fixed point of $K$. Hence all 96 exceptional
lines on $\tilde X$ are in the fixed point locus of $K$. There are exactly 12
orbits under the action of the group $K$.
Now  we have to count only the one dimensional fixed curves
in $\cal X$. This can be done with the results of \cite{FS2}.
 We just give the result: There are 12 elements of $K$ having
a one dimensional fixed point locus and each of them has 4 components,
which are elliptic curves. These are in the two K-orbits.
 \begin{lm}\label{l:fixedpts}
The number of components of the fixed point locus of $K$ on
$\tilde{\mathcal X}/K$ is $36$. 
 \end{lm}
\noindent
Now we are able to compute the Picard number of a Calabi-Yau model
of $\cal X/K$. The Picard number of the regular locus can be
computed by means of the results of section 6, especially
theorem 6.4 in \cite{FS2}. The result of a computation is 4. Hence we get:
 \begin{lm}
 The Picard number of a Calabi-Yau model of $\cal Y$ 
  is $40$. 
  \end{lm}
\noindent  
Let us   compute  the Euler number.  We recall that the  crepant resolution $\tilde{\cal X}$  has Euler number equal to $ 64$.  Since $K$ is abelian, the string
theoretic formula gives
$$ 
e(\tilde{\cal Y} )=\frac{1}{32} \sum_{(g,h)\in K\times K}e(\tilde{\cal X}^{<g,h>})=$$
$$\frac{64}{32} +\frac{3}{32} \sum_{g\neq id}e(\tilde{\cal X}^{g})
+\frac{1}{32} \sum_{\id\neq g\neq h\neq id}e(\tilde{\cal X}^{<g,h>})$$
Since the fixed point set of a single involution is
an elliptic curve or one of the 96 exceptional lines, we get
$$e=20+\frac{1}{32}\sum_{\id\neq g\neq h\neq id}e(\tilde{\cal X}^{<g,h>}).$$
We still have to discuss how for two different $g,h$,
which are different from the identity, the fixed point loci
intersect in $\tilde{\cal X}$. We want to compare this with the intersection
of the fixed point loci on the singular model $\cal X$.
We have to discuss two cases,  
\begin{itemize}
\item $g$ fixes a curve in $\cal X$ and $h$
fixes a node.
 \item\it  Both $g$ and $h$ have one dimensional fixed
point locus on $\cal X$ (4 elliptic curves). 
\end{itemize}
In the first case
there are 12 $g$ which fix a curve and 6 $h$ which fix a node.
Hence we have 72 cases to consider. In 48 cases the intersection
of the fixed point loci in $\cal X$ is empty.
Hence only 24 pairs are of interest. In each case
the fixed locus $\Fix(g)$ of $g$ is the union of 4
smooth elliptic curves
$$\Fix(g)=E_1\cup E_2\cup E_3\cup E_4.$$
and the fixed point locus of of $h$ consists of 16 nodes.
The intersection of $\Fix(g)$ and $\Fix(h)$ consists of 8 nodes.
Each single $E_i$ contains 4 of these 8 nodes.
This shows that in each of the 8 nodes two of the 4 elliptic
curves come together. Now we consider $\tilde{\cal X}$.
Since the fixed point set of $g$ is smooth, it consists of four
elliptic curves $\tilde E_1,\dots,\tilde E_4$, such that
the natural projection $\tilde E_i\to E_i$ is biholomorphic.
Let $a$ be one of the 8 nodes in $\Fix(g)\cap \Fix(h)$.
We can assume that $E_1,E_2$ are the two elliptic curves which
run into $a$.
Let $C$ be the
exceptional line over $a$.
Then $g$ induces an automorphism
of $C$ of order two. Since an involution $P^1$ has two fixed
points, we see that $\tilde E_1$ and $\tilde E_2$
each hit $C$ in one intersection point and both points are different.
So each of the 8 exceptional lines carries two
intersection points. This shows:
\begin{lm}
{Let $g\in K$ be an element with a one dimensional fixed point set,
and $h\in K$ an element, which fixes nodes. There are $24$
possibilities. The joint fixed point locus on $\tilde{\cal X}$
consists of $16$ points.}
\end{lm}
\noindent
In the formula for the Euler number each pair $(g,h)$ of the above
form contributes with $16/32$. We have 24 pairs. Together with the
pairs $(h,g)$ we get the contribution $24$ to the Euler number.
Hence we have
$$e=44+ \frac {1}{ 32}\sum_{\id\neq g\neq h\neq \id\atop
dim \Fix(g)=dim \Fix(h)=1}e(\tilde{\cal X}^{<g,h>}).$$ 
In the  second case, both $g$ and $h$ have one dimensional fixed
point locus on $\cal X$ (4 elliptic curves).
The number of intersection points of $\Fix(g)$ and $\Fix(h)$ on
$\cal X$ is 0, 8 or 16. The number of pairs $(g,h)$ with 8 intersection points
is 24 and that with 16 intersection points is 48.
\smallskip

Let us consider  pairs with $16$ intersection points. 
In this case one can check that none of the 16 is a node, and one can
check furthermore that   the contribution to the Euler
for each such pair is $(1/32)\cdot 16=1/2$.

\smallskip

Now we consider  pairs with $8$ intersection points. 
In this case one can check that all 8 intersection points are nodes.
Let $a$ be such a node. One can see that that two of the components of
$\Fix(g)$ run into $a$ and the same is true for $\Fix(h)$.
Moreover a simple computation gives that
$gh$ has $a$ as isolated singularity.
Hence as in the first case above $g$ has two fixed points $a_1,a_2$
on the exceptional line $C$ over $a$ and $h$ has the same    fixed
points.  
Hence 12 is the contribution to the Euler number.
We get as  contribution
$36=24+12$ to the Euler number. This gives
$$e=80$$
for the Euler number.
\begin{thm}\label{thm:hodge}
 A Calabi-Yau model of $\mathcal Y$
has Hodge numbers $h^{11}=40$, $h^{12}=0$.
\end{thm}

\noindent

\section{Explicit Calabi--Yau model of $\tilde{\mathcal Y}$}
\noindent
In this section we shall give alternative description of the
Calabi--Yau manifold $\mathcal {\tilde Y}$ using only the
equations~\eqref{eq:equations} 
of  $\mathcal Y$  as a complete intersection in the weighted
projective space $\Pj(1,1,1,1,2,2)$. These equations allows us to
consider $\mathcal Y$ as a $\Z/2\Z\oplus\Z/2\Z$ covering of the
projective space $\Pj^3$ branched along a pair of quartic surfaces. 
As a consequence we are able to use the standard methods of double
coverings to describe a crepant resolution of $\mathcal Y$, compute
its Euler characteristic and Hodge numbers (via the dimension of the
deformations space). We also give an explicit correspondences with the
van Geemen's and Nygaard's variety and the self fiber product of 
Beauville's surface.

Subtracting twice the first equation in \eqref{eq:equations} 
from the second one and changing the coordinate system 
$$(x_{0},x_{1},x_{2},x_{3},y_{4},y_{5})\mapsto(x_{0},x_{1},x_{2},x_{3},\tfrac 12(y_{4}+x_{2}^{2}-x_{0}^{2}-x_{1}^{2}-x_{3}^{2}))$$ we get the
following representation of $\mathcal Y$ as a complete intersection in
$\mathbb P(1,1,1,1,2,2)$
\[\def\arraycolsep{.3mm}
\begin{array}{ccl}
  y_{5}^{2}&=&x_{0}x_{1}x_{2}x_{3}\\
  y_{4}^{2}&=&(x_{0}+x_{1}+x_{2}+x_{3})\times
(x_{0}-x_{1}-x_{2}+x_{3})\times\\&&\times(x_{0}-x_{1}+x_{2}-x_{3})
\times(x_{0}+x_{1}-x_{2}-x_{3})
\end{array}
\]
Description of the rings of modular forms for varieties $\mathcal X$
and $\mathcal Y$ yields the following quotient map
\[
(X_{0},X_{1},X_{2},X_{3},Y_{0},Y_{1},Y_{2},Y_{3})\mapsto
(Y_{0}^{2},Y_{1}^{2},Y_{2}^{2},Y_{3}^{2},16X_{0}X_{1}X_{2}X_{3},Y_{0}Y_{1}Y_{2}Y_{3})
\]
so the action on $\mathcal X$ is diagonal given by the following group
\[K:=\{\varepsilon\in(\Z/2\Z)^8:\varepsilon_{0}=1,
\varepsilon_{1}\varepsilon_{2}\varepsilon_{3}=1,
\varepsilon_{4}\varepsilon_{5}\varepsilon_{6}\varepsilon_{7}=1\}\cong(\Z/2\Z)^5.\]

\noindent We are going to describe an explicit crepant resolution of
$\mathcal Y$.
Variety $\mathcal Y$ may be considered as $\Z
_{2}\oplus\Z _{2}$ covering of 
$\mathbb P^{3}$, or as an iterated double covering. The branch locus
consists of two quartics 
\begin{eqnarray*}
  D_{1}&=&\{x_{0}x_{1}x_{2}x_{3}=0\},\\ 
D_{2}&=&\{(x_{0}+x_{1}+x_{2}+x_{3})\times
(x_{0}-x_{1}-x_{2}+x_{3})\times\\&&\times(x_{0}-x_{1}+x_{2}-x_{3})
\times(x_{0}+x_{1}-x_{2}-x_{3})=0\}.
\end{eqnarray*}
 Both quartics $D_{1}$ and $D_{2}$ are sums of four faces
of tetrahedra 
in $\mathbb P^{3}$, so each of them gives four triple point and six 
double lines which we denote $l^{(1)}_{1},\dots,l^{(1)}_{6}$ and 
$l^{(2)}_{1},\dots,l^{(2)}_{6}$.

 Each of the lines $l^{(1)}_{i}$
intersect two of the lines $l^{(2)}_{j}$ giving rise to 12 fourfold
points of the octic $D:=D_{1}+D_{2}$, which we denote
$P_{1},\dots,P_{12}$.

The intersection $D_{1}\cap D_{2}$ is a sum of
sixteen lines (intersections of pair of planes a component of $D_{1}$
and a component of $D_{2}$) 
\[D_{1}\cap D_{2}=\sum_{i=1}^{16}C_{i}.\]
Let $\sigma_{1}:\widetilde{\mathbb P^{3}}\longrightarrow\mathbb P^{3}$ be the blow--up of $\mathbb P^{3}$ at
points $P_{1},\dots,P_{12}$, let $\tilde l^{(i)}_{j}$ denotes the strict
transform of $l^{(i)}_{j}$ and $\tilde D_{i}$ the strict transform of
$D_{i}$. Then the lines $\tilde l^{(1)}_{i}$ and $\tilde l^{(2)}_{j}$ are disjoint
whereas any three out of $\tilde l^{(1)}_{i}$  and any three out of
$\tilde l^{(2)}_{i}$     intersect at a triple point. Moreover we have
$\tilde D_{i}=\sigma_{1}^{*}D_{i}-2\operatorname{exc}(\sigma_{1})$,
$K_{\widetilde{\mathbb P^{3}}}=\sigma_{1}^{*}K_{\mathbb P^{3}}+2\operatorname{exc}(\sigma_{1})$
hence
\[K_{\widetilde{\mathbb P^{3}}}+\tfrac12 (\tilde D_{1}+\tilde
D_{2})=\sigma_{1}^{*}(K_{\mathbb P^{3}}+\tfrac12(D_{1}+D_{2})) 
.\]
Let $\sigma_{2}:\mathbb P^{*}:\longrightarrow \widetilde{\mathbb P^{3}}$ be the composition of
blow--ups of (strict transforms of) lines $\tilde l^{(i)}_{j}$. 
For each blow--up the strict transform of the quartic which contain it
equals the pullback minus twice the exceptional divisor, whereas for
the other quartic the strict transform equals the pullback.

Denote
by
\[\sigma:\mathbb P^{*}\longrightarrow \mathbb P^{3}\]
composition $\sigma:=\sigma_{2}\circ\sigma_{1}$ and by $D_{i}^{*}$ the
strict transform of $D_{i}$. Then $D_{1}^{*}$ and $D_{2}^{*}$ are smooth
divisors intersecting transversally along a disjoint sum of 16 lines. 

Let $\tilde {\mathcal Y}$ be a $\Z/2\Z\oplus\Z/2\Z$
Galois covering 
of $\mathbb P^{*}$ branched along divisors $D_{1}^{*}$ and $D_{2}^{*}$.
\begin{lm}\label{l:bi-double}
    \begin{eqnarray*}
\pi_{*}\mathcal O_{\tilde {\mathcal Y}}&=&\mathcal
O_{\mathbb P^{*}}(-\tfrac12(D_{1}^{*}+D_{2}^{*}))\oplus\mathcal
O_{\mathbb P^{*}}(-\tfrac12D_{2}^{*})
  \oplus\mathcal O_{\mathbb P^{*}}(-\tfrac12D_{1}^{*})
  \oplus\mathcal O_{\mathbb P^{*}}.   \\
  \pi_{*}\Theta_{\tilde{\mathcal Y}}&=&\Theta_{\mathbb P^{*}}
  (-\tfrac12D^{*})\oplus\Theta_{\mathbb P^{*}}(\log
  D_{1})(-\tfrac12D_{2}^{*})\oplus\\
  &&\oplus\Theta_{\mathbb P^{*}}(\log D_{2}^{*})(-\tfrac12D_{1}^{*})
  \oplus\Theta_{\mathbb P^{*}}(\log D^{*}).\\
  K_{\tilde{\mathcal Y}}&=&0.
  \end{eqnarray*}

\end{lm}
\begin{proof}
  The first two assertion can be directly verified in local
  coordinates, they also follows from factoring the map $\pi$ into a
  composition of two double covering: double covering of $\mathbb P^{*}$
  branched along $D_{1}^{*}$ followed by a double covering branched
  along pullback of $D_{2}^{*}$ (or similar with $D_{1}$ and $D_{2}$
  exchanged). From this factorization it follows that
  $K_{\tilde{\mathcal Y}}=K_{\mathbb P^{*}}+\tfrac12
  (D_{1}^{*}+D_{2}^{*})=\pi^{*}(K_{\mathbb P^{3}}+\tfrac12(D_{1}+D_{2}))=0$.
 
\end{proof}

Now, we can give another proof of Thm.~\ref{thm:cymodel} and
Thm.~\ref{thm:hodge}.
Since the map $\sigma$ is a composition of blow--ups with smooth
centers $\sigma_{*}\mathcal O_{\mathbb P^{*}}=\mathcal O_{\mathbb P^{3}}$ and
$R^{i}\sigma_{*}\mathcal O_{\mathbb P^{*}}=0$, for $i>0$. So by the Leray spectral
sequence and Serre duality  $H^{1}(\mathcal O_{\mathbb P^{*}})=H^{1}(\mathcal O_{\mathbb P^{3}})=0$
and
$H^{1}(\mathcal O_{\mathbb P^{*}}(-\tfrac12(D_{1}^{*}+D_{2}^{*})))=
H^{1}(K_{\mathbb P^{*}})=H^{2}(\mathcal O_{\mathbb P^{*}})=H^{2}(\mathcal O_{\mathbb P^{3}})=0$.   

\smallskip\noindent
{\bf Claim.} $\sigma_{*}\mathcal 
O_{\mathbb P^{*}}(-\tfrac12D_{2}^{*})=\mathcal O_{\mathbb P^{3}}(-\tfrac12D_{2})$,  
$R^{i}\sigma_{*}\mathcal  O_{\mathbb P^{*}}(-\tfrac12D_{2}^{*})=0$ for
$i>0$.  
To prove the claim we shall consider every blow--up separately, let
$\mathcal L$ be a line bundle on a smooth threefold $P$ and let 
$\tau :\tilde P\longrightarrow P$ be a blow--up of a smooth subvariety $C\subset
P$ with exceptional divisor. Let $M$ be a line bundle on $\tilde P$
satisfying one of the following three conditions
\begin{itemize}
\item $C$ is a curve and $\mathcal M=\tau^{*}\mathcal L$,
\item $C$ is a curve and $\mathcal M=\tau^{*}\mathcal L\otimes
  \mathcal O_{\tilde P}(E)$,
\item $C$ is a point and $\mathcal M=\tau^{*}\mathcal L\otimes
  \mathcal O_{\tilde P}(E)$.
\end{itemize}
 In the first case, by the projection formula, $\tau _{*} \mathcal
 M=\mathcal L$ and $R^{i}\tau_{*}\mathcal M=0$. 
In the other two cases consider the following exact sequence
\[0\longrightarrow\tau^{*}\mathcal L\longrightarrow\mathcal M\longrightarrow\tau^{*}\mathcal L\otimes
\mathcal O_{E}(-1) \longrightarrow 0.
\]
Since $\tau_{*}(\mathcal O_{E}(-1))=R^{i}\tau_{*}(\mathcal O_{E}(-1))=0$, applying
the direct image functor to the above exact sequence yields
$\tau _{*} \mathcal
 M=\mathcal L$ and $R^{i}\tau_{*}\mathcal M=0$ and the claim follows.

From the Leray spectral sequence we get 
$$H^{1}(\mathcal 
O_{\mathbb P^{*}}(-\tfrac12D_{2}^{*}))=H^{1}(\sigma_{*}\mathcal 
O_{\mathbb P^{*}}(-\tfrac12D_{2}^{*}))=H^{1}(\mathcal O_{\mathbb P^{3}}(-\tfrac12D_{2}))=0$$
and (by symmetry) $H^{1}(\mathcal 
O_{\mathbb P^{*}}(-\tfrac12D_{2}^{*}))=0$.

The map $\pi$ is finite so using Lemma~\ref{l:bi-double} we get 
\[H^{1}(\mathcal O_{\tilde{\mathcal Y}})=0\]
which proves that $\tilde{\mathcal Y}$ is a Calabi--Yau threefold.

By the above description $\mathbb P^{*}$ is the projective space $\mathbb P^{3}$
blown--up at twelve points and twelve lines so
 \[e(\mathbb P^{*})=4+12\times 2+12\times 2=52.\]
Observe that  blowing--up a double
line containing  a triple point blows--up also one of the planes
containing this point, whereas blowing--up a fourfold point blows--up
all four planes through theis point.  Consequently $D_{1}^{*}$ is a sum of four
planes blown--up 28 times, so 
\[e(D_{1}^{*})=e(D_{2}^{*})=4\times 3+28=40\]
and $D_{1}^{*}\cap D_{2}^{*}$ is a disjoint sum of 16 lines so 
\[e(D_{1}^{*}\cap D_{2}^{*})=32.\]
Now,
\[e(\tilde {\mathcal
  Y})=4e(\mathbb P^{*})-2e(D_{1}^{*})-2e(D_{2}^{*})+e(D_{1}^{*}\cap
D_{2}^{*}) =4\times52-2\times 80+32=80.
\]
To prove that $h^{1,2}(\tilde{\mathcal Y})=0$, we shall proceed as in
\cite{CvS}. By \cite[Thm.~4.7]{CvS} $H^{1}(\Theta_{\mathbb P^{*}}(\log
D^{*}))$ is isomorphic to the space of equisingular deformations of
$D$ in $\mathbb P^{3}$, moreover it is isomorphic to $(I_{eq}(D)/J_{F})_{8}$,
where $J_{F}$ is the jacobian ideal of $D$ and
\[I_{eq}=\bigcap_{i=1}^{12}(I(P_{i})^{4}+J_{F})\cap\bigcap
_{i=1}^{6}\bigcap_{j=1}^{2}(I(l_{j}^{(i)})^{2}+J_{F}) 
\]
is the equisingular ideal. Using this formula we check with Singular
(\cite{Singular}) that $H^{1}(\Theta_{\mathbb P^{*}}(\log
D^{*}))=0$. 

As in the resolution of $\mathcal Y$ we blow--up only rational curves,
by \cite[Prop.~5.1]{CvS} $H^{1}(\Theta_{\mathbb P^{*}}(-\tfrac12D^{*}))=0$.

Consider the following exact sequence
\[
0\longrightarrow\Theta_{\mathbb P^{*}}(\log
D_{1}^{*})(-\tfrac12D_{2}^{*})\longrightarrow\Theta_{\mathbb P^{*}}(-\tfrac12D_{2}^{*}) 
\longrightarrow\mathcal N _{D_{1}^{*}}(-\tfrac12D_{2}^{*})\longrightarrow0.
\]
We shall study first $\Theta_{\mathbb P^{*}}(-\tfrac12D_{2}^{*})$ and again
consider separately a single blow--up $\tau:\tilde P\longrightarrow P$ with a
smooth center $C$. We have the same three cases 
\begin{itemize}
\item $C$ is a curve and $\mathcal M=\tau^{*}\mathcal L$,
\item $C$ is a curve and $\mathcal M=\tau^{*}\mathcal L\otimes
  \mathcal O_{\tilde P}(E)$,
\item $C$ is a point and $\mathcal M=\tau^{*}\mathcal L\otimes
  \mathcal O_{\tilde P}(E)$,
\end{itemize}
where $\tau$ is as before, and consider the vector bundle
$\Theta_{\tilde P}\otimes \mathcal M$. Using \cite[Sect.~5]{CvS} in
the first and third cases ($k>0$ in notations of \cite{CvS}) we get 
$\tau_{*}(\Theta_{\tilde P}\otimes \mathcal M)=\Theta_{P}\otimes
\mathcal L$ and $R^{i}\tau_{*}(\Theta_{\tilde P}\otimes \mathcal
M)=0$. Since in this case $\mathcal N _{C}\otimes \mathcal L=K_{C}$, we get 
\[H^{1}(\Theta_{\mathbb P^{*}}(-\tfrac12D_{2}^{*}))=0.\]
Finally, to find $H^{0}(\mathcal N _{D_{1}^{*}}(-\tfrac12D_{2}^{*}))$ consider
the exact sequence
\[0\longrightarrow\mathcal O_{\mathbb P^{*}}(-\tfrac12D_{2}^{*})\longrightarrow
\mathcal O_{\mathbb P^{*}}(D_{1}^{*}-\tfrac12D_{2}^{*})\longrightarrow
\mathcal N _{D_{1}^{*}}(-\tfrac12D_{2}^{*})\longrightarrow0.\] 
Since $\sigma_{*}(\mathcal
O_{\mathbb P^{*}}(D_{1}^{*}-\tfrac12D_{2}^{*}))=\mathcal
O_{\mathbb P^{3}}(3)\otimes\mathcal I$, where $\mathcal I$ is the ideal of
functions vanishing at $P_{1},\dots,P_{12}$ and vanishing to order two
along $l^{(1)}_{1},\dots,l^{(1)}_{6}$, we get $H^{0}(\mathcal
O_{\mathbb P^{*}}(D_{1}^{*}-\tfrac12D_{2}^{*}))=0$. Since $H^{1}(\mathcal
O_{\mathbb P^{*}}(-\tfrac12D_{2}^{*}))=0$, we get 
$H^{0}(\mathcal N _{D_{1}^{*}}(-\tfrac12D_{2}^{*}))=0$ and consequently 
$H^{1}(\Theta_{\mathbb P^{*}}(-\tfrac12D^{*}))=0$. By the above exact
sequence we get $H^{1}(\Theta_{\mathbb P^{*}}(\log
D_{1}^{*})(-\tfrac12D_{2}^{*}))=0$ and (by symmetry) 
$H^{1}(\Theta_{\mathbb P^{*}}(\log
D_{2}^{*})(-\tfrac12D_{1}^{*}))=0$.

Since the map $\pi$ is finite Lemma~\ref{l:bi-double} yields
\begin{eqnarray*}
    H^{1}(\Theta_{\tilde{\mathcal Y}})=
    H^{1}(\Theta_{\mathbb P^{*}}
  (-\tfrac12D^{*}))+H^{1}(\Theta_{\mathbb P^{*}}(\log
  D_{1})(-\tfrac12D_{2}^{*}))+\\
  +H^{1}(\Theta_{\mathbb P^{*}}(\log D_{2}^{*})(-\tfrac12D_{1}^{*}))
  +H^{1}(\Theta_{\mathbb P^{*}}(\log D^{*}))=0
\end{eqnarray*}
and by the Serre duality
\[h^{1,2}(\tilde{\mathcal Y})=0.\]
Since the Hodge numbers of a Calabi--Yau manifold $\tilde{\mathcal Y}$
satisfy $e(\tilde{\mathcal Y})=2(h^{1,1}(\tilde{\mathcal
  Y})-h^{1,2}(\tilde{\mathcal Y}))$ we conclude 
\[h^{1,1}(\tilde{\mathcal Y})=40.\] 

\noindent There is another intersection of four quadrics related to the
Calabi--Yau manifold $\tilde{\mathcal Y}$. After the coordinate change 
\[(x_{0}:x_{1}:x_{2}:x_{3}:y_{4}:y_{5})\mapsto
(x_{0}+x_{1},x_{0}-x_{1},x_{2}+x_{3},x_{2}-x_{3}:y_{4}:\tfrac12y_{5})\] 
the equations are transformed into more symmetric
\begin{eqnarray*}
y_{5}^{2}&=&(x_{0}^{2}-x_{1}^{2})(x_{2}^{2}-x_{3}^{2}),\\
y_{4}^{2}&=&(x_{0}^{2}-x_{2}^{2})(x_{1}^{2}-x_{3}^{2}).  
\end{eqnarray*}
Consequently it is a $\Z/2\Z\oplus\Z/2\Z$--quotient of the following
intersection of four quadrics
\begin{equation}\label{verr}
  \begin{aligned}
u_{0}^{2}&=&x_{0}^{2}-x_{1}^{2},\\
u_{1}^{2}&=&x_{1}^{2}-x_{2}^{2},\\  
u_{2}^{2}&=&x_{2}^{2}-x_{3}^{2},\\
u_{3}^{2}&=&x_{3}^{2}-x_{0}^{2}    
  \end{aligned}
\end{equation}
in $\mathbb P^{7}$.
The intersection $S$ of two quadric in $\mathbb P^{4}$
\begin{eqnarray*}
u_{0}^{2}&=&x_{0}^{2}-x_{1}^{2},\\
u_{1}^{2}&=&x_{1}^{2}-x_{2}^{2}\\  
\end{eqnarray*}
is singular at points $(0:0:1:0:\pm i)$, $(1:0:0:\pm1:0)$, the
rational map $\pi:S\ni(x_{0}:x_{1}:x_{2}:u_{0}:u_{1})\longrightarrow(x_{0}:x_{2})$ is
undetermined at points $(0:1:0:\pm i:\pm1)$
(intersection of the surface $S$ with the plane
$x_{0}=x_{2}=0$). Blowing--up
$S$ at singular points and then at points of indeterminacy  yields a
rational elliptic surfaces $\tilde\pi:\tilde S\longrightarrow\mathbb P^{1}$  with
fours singular fibers: of type $I_{4}$ at $0,\infty$ and $I_{2}$ at
$\pm1$. It means that $\tilde S$ is the Beauville modular surfaces
associated to the group $\Gamma_{1}(4)\cap \Gamma(2)$ and the
intersection~\eqref{verr} is the self fiber product of $\tilde S$. \smallskip
 
\smallskip
From the above description it follows that  $\tilde{\mathcal
  Y}$  is modular with the unique cusp form of weight 4 and level
8. One can  also prove that using the Faltings--Serre--Livn\'e
method. Using a computer program we verify that for $p$ prime,
$p\le97$ the number of points in $\mathcal X(\mathbb F_{p})$ equals 
\[1+p^{3}-a_{p}+16(p+p^{2})-12(2p+p^{2}),\]
where $a_{p}$ is the coefficient of the cusp form.

\bigskip

\section{K3 fibration and the Picard group}

\noindent  The Hodge number $h^{1,1}(\tilde{\mathcal Y})=40$ equals the
Picard number of the Calabi--Yau manifold $\tilde{\mathcal Y}$.
The resolution of singularities of $\tilde{\mathcal Y}$ yields 37
apparent linearly independent divisors:

\begin{itemize}
\item pullback of a hyperplane section in $\mathbb P^{3}$,
\item 12 blow--ups $\Z/2\Z\oplus\Z/2\Z$ covers of a plane the
  exceptional loci of blow--ups of fourfold points,
\item 24 blow--ups of double covers of exceptional divisors of
  blow--up of a double line, since after blowing--up fourfold points
  any double line is disjoint from one of the branch divisors, the
  $\Z/2\Z\oplus \Z/2\Z$ covers splits into a pair of double covers.
\end{itemize}

\begin{rem}The twelve fourfold points points have the form
  $(1:\pm1:0:0:0:0)$ and their permutations of $x_{1},\dots,x_{4}$ coordinates. 
By the description of the quotient map they  correspond to the  12
orbits of the nodes under $K$ action.
 
The twelve lines $l^{(i)}_{j}$ corresponds by the quotient map to the
intersections of $\mathcal X$ with linear subspaces $X_{k}=X_{l}=0$
or $Y_{k}=Y_{l}=0$ which are sums of four elliptic curves.  

So the above description of 36 linearly independent divisors agrees
with the description given in Lemma~\ref{l:fixedpts}.
\end{rem}
%
\noindent
In this way we can identify rank 37 subgroup in the Picard group. To
identify the remaining divisors we can use one of the K3 fibrations
on  $\tilde{\mathcal Y}$. Fix a double line of one of the quartics
(f.i. fix the line
$m:=\{(x_{0}:x_{1}:x_{2}:x_{3})\in\mathbb P^{3}:x_{2}=x_{3}=0\}\subset D_{1}$)
and let
$P_{(s:t)}:=\{(x_{0}:x_{1}:x_{2}:x_{3})\in\mathbb P^{3}:sx_{2}+tx_{3}=0\}$
($(s:t)\in\mathbb P^{1}$) be the pencil of planes  that defines a fibration
on $\tilde{\mathcal Y}$. For $(s:t)\not=0,\infty,\pm1$, the fiber
$S_{(s:t)}$ is a 
smooth K3 surface, it can be described as resolution of the complete
intersection in 
$\mathbb P(1,1,1,2,1)$
\[\def\arraycolsep{.3mm}
\begin{array}{ccl}
  y_{5}^{2}&=&x_{0}x_{1},\\
  y_{4}^{2}&=&(tx_{0}+tx_{1}+(t-s)x_{2})\times
(tx_{0}-tx_{1}+(-t-s)x_{2})\times\\&&\times(tx_{0}-tx_{1}+(t+s)x_{2})
\times(tx_{0}+tx_{1}+(-t+s)x_{2}).
\end{array}
\]
This is a $\Z/2\Z\oplus\Z/2\Z$ covering of $\mathbb
P^{2}$ branched along a pair and a quadruple of lines in general position,
the branch curves have seven nodes. Described resolution of singularities 
of $\mathcal Y$ induces also a resolution of singularities of a
generic fiber by blowing--up the double points of the branch curves.
Each of them induces two independent divisors in the Picard group,
together with a hyperplane section we get 16 linearly independent
divisors. 

There are however three more independent 
divisors, the lines $tx_{0}+sx_{1}=0$ and $sx_{0}+tx_{1}=0$ and the
conic $tx_{0}x_{1}+sx_{2}^{2}=0$ in $P_{(s:t)}$ intersects the branch
divisors only with multiplicity two, so they split in the covering
into four components. Taking one components from each of them shows
that the Picard number of the generic fiber is at least 19, which is
the biggest possible.

The singular fibers are reducible, comparing with the resolution we
get:
\begin{itemize}\leftskip=-7mm
\item [] the fiber $S_{{1:1}}$  (resp. $S_{{1:-1}}$) has three 3
components: the strict 
transform of the plane, divisor corresponding to the blow--up of the
point $(0:0:1:-1)$ (resp. $(0:0:1:1)$) and the line
$x_{0}+x_{1}=x_{2}+x_{3}=0$ (resp. $x_{0}+x_{1}=x_{2}+x_{3}=0$),\medskip
\item [] the fiber $S_{(1:0)}$ (resp. $S_{(0:1)}$) has 9 components:
the strict transform of the plane, four divisor corresponding to the
blow--up of points $(0:1:0:1)$, $(0:1:0:-1)$, $(1:0:0:1)$,
$(1:0:0:-1)$ (resp. $(0:1:1:0)$, $(0:1:-1:0)$, $(1:0:1:0)$,
$(1:0:-1:0)$) four divisors (two pairs) corresponding to the lines 
$x_{0}=x_{2}=0$ and $x_{1}=x_{2}=0$ (resp. $x_{0}=x_{2}=0$ and
$x_{1}=x_{2}=0$). \smallskip
\end{itemize}

\noindent On the Calabi--Yau model the three divisors on the generic fiber of
fibration correspondfs to components of the strict transforms of the
quadrics
\[x_{0}x_{1}=x_{2}x_{3}, x_{0}x_{2}=x_{1}x_{3}, 
x_{0}x_{3}=x_{1}x_{2}\]
in  $\mathbb P^{3}$.
\begin{rem}
  Since 
\[Y_{0}^{2}Y_{1}^{2}-Y_{2}^{2}Y_{3}^{2}=4(X_{0}X_{2}-X_{1}X_{3})^{2},\]
components of the strict transform of the quadric
$x_{0}x_{1}=x_{2}x_{3}$ correspond via the quotient map to the
components of the intersection of 
  $\mathcal X$ with  the quadric $X_{0}X_{2}-X_{1}X_{3}$. These Weil
  divisors on $\mathcal X$ are not $\mathbb Q$--Cartier, they give a
  projective small resolution of van Geemen's and Nygaard's variety
  (cf. \cite{FS2}). 
\end{rem}

 \end{document}